\documentclass[a4paper]{amsart}
\usepackage[latin1]{inputenc}
\usepackage[english]{babel}
\usepackage{amssymb}

\makeatletter

 \theoremstyle{plain}
 \newtheorem{thm}{Theorem}
 \newtheorem{prop}[thm]{Proposition}
 \newtheorem{lem}[thm]{Lemma}
 
 \newtheorem{cor}[thm]{Corollary}
 \newtheorem{rem}[thm]{Remark}

\usepackage{amssymb}

\newcommand{\Vertex}{\operatorname{Vert}}
\newcommand{\Edge}{\operatorname{Edge}}
\newcommand{\val}{\operatorname{val}}
\newcommand{\rk}{\operatorname{rk}}
\newcommand{\ev}{\operatorname{ev}}
\newcommand{\pd}{\operatorname{P.D.}}
\newcommand{\pt}{\operatorname{pt}}

\newcommand{\DM}[1][m]{\overline{\mathcal{M}}_{0,#1}}
\newcommand{\mfix}{\mathcal{M}_\Gamma}
\newcommand{\bsp}{\mathbb{P}_{\mathbb{P}^1}(\mathcal{O}^{\oplus (r-2)}\oplus\mathcal{O}(1)\oplus\mathcal{O}(1))}
\makeatother

\begin{document}

\title{Multiple quantum products in toric varieties}

\author{Holger Spielberg}

\date{\today}

\address{Centro de Analíse Matemática, Geometria e Sistemas
Dinâmicos,
Instituto Superior Técnico,
Av. Rovisco Pais,
1049-001 Lisboa,
Portugal}

\email{Spielberg@member.ams.org}

\begin{abstract}
We generalize the author's formula for Gromov--Witten invariants of symplectic
toric manifolds (see \cite{spielberg-phd}, \cite{spielberg-GWcras}) to those
needed to compute the quantum product of more than two classes directly, i.e.
involving the pull--back of the Poincar\'e dual of the point class in the
Deligne--Mumford spaces $ \DM $.
\end{abstract}
\maketitle

\section{Introduction}

Let $ (X,\omega ) $ be a symplectic manifold with compatible almost--complex
structure $ J $. If $ g $ and $ m $ are non--negative integers, we
denote by $ \overline{\mathcal{M}}_{g,m} $ the Deligne--Mumford space of
genus--$ g $ curves with $ m $ marked points. If furthermore $ A\in H_{2}(X,\mathbb Z) $
denotes a degree--$ 2 $ homology class of $ X, $ $ \mathcal{M}_{g,m}(X,A) $
will be the space of stable genus--$ g $ $ J $--holomorphic maps to $ X $
with homology class $ A $. The Gromov--Witten invariants of $ X $ (see
for example \cite{ruan-tian-mathematicaltheory}, \cite{behrend-gwinvariants})
are multi--linear maps
\[
\Phi _{g,m}^{A,X}:
H^{*}(\overline{\mathcal{M}}_{g,m},\mathbb Q)\otimes H^{*}(X,\mathbb Q)^{\otimes m}
\longrightarrow \mathbb Q
\]
that are defined as follows. Let
$ \pi :\mathcal{M}_{g,m}(X,A)\longrightarrow \overline{\mathcal{M}}_{g,m} $
be the natural projection map forgetting the map to X (and stabilizing the curve
if necessary). Furthermore, let $ \ev_{i}:\mathcal{M}_{g,m}(X,A)\longrightarrow X $
be the evaluation map at the $ i^{th} $ marked point, that is the map that
sends a stable map $ (C;x_{1},\ldots ,x_{m};f) $ to $ f(x_{i})\in X $.
Then for classes $ \beta \in H^{*}(\overline{\mathcal{M}}_{g,m},\mathbb Q) $
and $ \alpha _{i}\in H^{*}(X,\mathbb Q) $, the Gromov--Witten invariants
are defined by
\begin{equation}
\label{eq:gwdef}
\Phi _{g,m}^{A,X}(\beta ;\alpha _{1},\ldots ,\alpha _{m})=
\int _{[\mathcal{M}_{g,m}(X,A)]^{\text {virt}}}
\pi ^{*}(\beta )\wedge \ev_{1}(\alpha _{1})\wedge \ldots \wedge \ev_{m}(\alpha _{m}).
\end{equation}
Here the integration on the right hand side is not over the entire moduli space
but over the so-called virtual fundamental class.

For the case of $ (X,J) $ being a smooth projective variety with a $ (\mathbb C^{*}) $--action,
Graber and Pandharipande have proven (see \cite{graber-pandharipande}) that
Bott--style localization techniques apply to the integral in (\ref{eq:gwdef}).
Their techniques can easily be extended to torus actions, so in particular they
apply to smooth projective toric varieties (for toric varieties see for example
\cite{fulton-toricvarieties}, \cite{oda-convexbodies}).

In \cite{spielberg-phd} (also see \cite{spielberg-GWcras}), using these localisation
techniques for the virtual fundamental class, we have proven an explicit combinatorial
formula of the genus--$ 0 $ Gromov--Witten invariants for smooth projective
toric varieties for the cases when $ \beta =1\in H^{0}(\DM,\mathbb Q) $.
In this note we will derive a similar formula for the case where the class $ \beta  $
is the maximal product of (the Chern class of) cotangent lines to the marked points,
that is for those classes $\beta$ which are Poincar\'e dual to a finite number of points
in $\DM$.

Knowing the Gromov--Witten invariants for $\beta=\pd(\pt)$
makes computations of products in the (small) quantum cohomology ring easier.
For two cohomology classes $ \gamma _{1},\gamma _{2}\in H^{*}(X,\mathbb Q) $,
their quantum product is defined to be
\[
\gamma _{1}\star \gamma _{2}=
\sum _{A\in H_{2}(X,\mathbb Z)}\sum _{i}
\Phi _{0,3}^{X,A}(\gamma _{1},\gamma _{2},\delta _{i})\delta _{i}^{\vee }q^{A}
\]
where the inner sum runs over a basis $ (\delta _{i}) $ of $ H^{*}(X,\mathbb Q) $.
Here $ (\delta _{i}^{\vee }) $ denotes the basis of $ H^{*}(X,\mathbb Q) $
dual to $ (\delta _{i}) $. It is easy to show that the product of more than
two classes is given by
\begin{equation}\label{eq:multquantumproduct}
\gamma _{1}\star \cdots \star \gamma _{r}=
\sum _{A\in H_{2}(X,\mathbb Z)}\sum _{i}
\Phi _{0,r+1}^{X,A}(\pd(\pt);\gamma _{1},\ldots ,\gamma _{r},\delta _{i})
\delta _{i}^{\vee }q^{A}
\end{equation}
 where $ \pd(\pt) $ denotes the Poincaré dual of a point in $ \DM[r+1] $.
By Witten's conjecture (see \cite{witten-2dimgravity}, proven by Kontsevich
in \cite{kontsevich-Airyfunction}), we know that (a multiple of) the class $ \pd(\pt) $
in $ H^{*}(\overline{M}_{0,r+1}) $ can be expressed as the product of cotangent
line classes, hence the invariants in (\ref{eq:multquantumproduct}) can be
computed directly by the formula proposed in this note.

The techniques used in this note also yield the invariants for $\beta$ being
a different product of Chern classes of such cotangent line bundles. However, it
seems to be much more to difficult to formalize such a more general approach.
For the sake of a (hopefully) better exposition of the key ideas, we
leave the more general case to the interested reader.

The structure of the paper is as follows: In Section 2 we recall some results
on toric varieties, mostly to fix our notation. In Section 3 we quickly describe
the fixed point components of $\mathcal{M}_{0,m}(X_\Sigma,A)$ with respect
to the action induced from $X_\Sigma$. In Section 4 we recall the localization
results for toric varieties, to apply them in Section 5 to the case where
$\beta = \pd(\pt)$. In Section 6 we finally give the formula for the
Gromov--Witten invariants for symplectic toric manifolds in this case, and in
Section 7 we illustrate the formula on the example of $\bsp$; as an interesting
byproduct, we derive the quantum cohomology ring of this variety (also using
recent results of \cite{costa-miro}), which surprisingly coincides with
Batyrev's ring stated in \cite{batyrev-quantumcohomologytoric}.

\section{Preliminaries of toric varieties}

We will quickly recall some facts about toric varieties and mostly introduce
our notation --- our standard references for this section are
\cite{batyrev-quantumcohomologytoric}, \cite{fulton-toricvarieties} and
\cite{oda-convexbodies}.

Let $ X_{\Sigma } $ be a smooth projective toric variety of complex dimension
$ d $, given by the fan $ \Sigma  $. Choose a class $ \omega  $ in
the Kähler cone of $ X_{\Sigma } $, and let $ \Delta _{\omega } $ be the
corresponding moment polytope. On the variety $ X_{\Sigma } $, the $ d $-dimensional
torus $ T_{N}:=(\mathbb C^{*})^{d} $ acts effectively, and the (irreducible)
subvarieties of $ X_{\Sigma } $ that are left invariant under this action
are in one--to--one correspondence with the facets of the polytope $ \Delta _{\omega } $.
Moreover, the $ T_{N} $--invariant divisors (which are in one--to--one correspondence
to faces of $ \Delta _{\omega } $) generate the cohomology ring $ H^{*}(X_{\Sigma },\mathbb Z) $
of $ X_{\Sigma } $ --- we will denote the faces of $ \Delta _{\omega } $
by $ Z_{1},\ldots ,Z_{n} $. We also remind the reader, that the relations
between these divisors in the cohomology ring are given by the combinatorics
of $ \Delta _{\omega } $ or equivalently, by that of the fan $ \Sigma  $.
For (higher--degree) cohomology classes we will sometimes use multi--index notation,
i.~e.~$ Z^{l} $ expands to $ Z_{1}^{l_{1}}\cdots Z_{n}^{l_{n}} $.

We will be using the weights of the torus action on the tangent bundle at fixed
points of $ X_{\Sigma } $. The vertices of $ \Delta _{\omega } $ are in
one--to--one correspondence with these fixed points, and we will usually denote
these vertices by the greek letter $ \sigma  $.
For any vertex $\sigma$, there are exactly $ d $
edges $ e_{1},\ldots ,e_{d} $ in $ \Delta _{\omega } $ that meet at $ \sigma  $.
Each edge of the polytope $ \Delta _{\omega } $ correspondents
to $ T_{N} $--invariant $ \mathbb C\mathbb P^{1} $ in $ X_{\Sigma } $.
Then the tangent space $ T_{\sigma }X_{\Sigma } $ at $ \sigma  $ $ T_{N} $-invariantly
splits into the tangent lines along these subvarieties. If we denote by
$ \sigma _{1},\ldots ,\sigma _{d} $
the vertices that are connected by the edges $ e_{1},\ldots ,e_{d} $ to $ \sigma  $,
we will denote by $ \omega _{\sigma _{i}}^{\sigma } $ the weight of the $ T_{N} $--action
on $ T_{\sigma }X_{\Sigma } $ into the direction of $ e_{i} $.

When referring to a degree--$ 2 $ homology class $ \lambda \in H_{2}(X_{\Sigma },\mathbb Z) $,
we will usually give its intersection vector $ (\lambda _{i})_{i=1,\ldots n} $
with the divisor classes $ Z_{i} $, that is $ \lambda _{i}:=\langle Z_{i},\lambda \rangle  $.
Note however, that the $ \lambda _{i} $ have to satisfy certain linear relations
to represent a degree--$ 2 $ homology class. In fact we have that
$ \dim H_{2}(X_{\Sigma },\mathbb Z)=\dim H^{2}(X_{\Sigma },\mathbb Z)=n-d $.

\section{Fixed--point components of the induced action on \protect$ \mathcal{M}_{0,m}(X_{\Sigma },A)\protect $}

Remember that an element of $ \mathcal{M}_{0,m}(X_{\Sigma },A) $ is (up to
isomorphisms) a tuple $ (C;x_{1},\ldots ,x_{m};f) $ where $ C $ is an
algebraic curve of genus zero with singularities at most ordinary double points,
$ x_{i}\in C_{\text {smooth}} $ are marked points, and $ f:C\longrightarrow X_{\Sigma } $
is the map to the variety $ X_{\Sigma } $. The $ T_{N} $--action on $ X_{\Sigma } $
then induces an action of $ T_{N} $ on $ \mathcal{M}_{0,m}(X_{\Sigma },A) $
by simple composition with the map $ f $, that is
$ t\cdot (C;x_{1},\ldots ,x_{m};f)=(C;x_{1},\ldots ,x_{m};t\circ f) $,
where
\[
t\circ f:C\stackrel{f}{\longrightarrow }X_{\Sigma }
\stackrel{\varphi _{t}}{\longrightarrow }X_{\Sigma }
\]
is the composition of $ f $ with the diffeomorphism $ \varphi _{t} $
given by the action of $ t $ on $ X_{\Sigma } $.

It is then easy to see (cf. \cite{spielberg-phd}), that the image of a fixed
point $ (C;\underline{x};f)\in \mathcal{M}_{0,m}(X_{\Sigma },A) $
must be left invariant by the $ T_{N} $--action, or in other words, it has
to live on the $ 1 $--skeleton of $ \Delta _{\omega } $. Moreover, the
marked points $ x_{i} $ of such a stable map will have to be mapped to
fixed points $ \sigma _{i} $ in $ X_{\Sigma } $. The fixed-point components
of $ \mathcal{M}_{0,m}(X_{\Sigma },A) $ can then be characterized by so--called
$ \mathcal{M}_{0,m}(X_{\Sigma },A) $-graphs $ \Gamma  $
(see \cite[Definition 6.4]{spielberg-phd})---these
are graphs on the $ 1 $--skeleton on $ \Delta _{\omega } $, without loops,
with decorations representing the position of the marked points and the multiplicities
of the map to $ X_{\Sigma } $ on the irreducible components. If $ \Gamma  $
is such a graph, we will usually denote by $ \mathcal{M}_{\Gamma } $ the
product of Deligne--Mumford spaces corresponding to the graph $ \Gamma  $,
which is up to a finite automorphism group isomorphic to the fixed--point component
in $ \mathcal{M}_{0,m}(X_{\Sigma },A) $ corresponding to $ \Gamma  $.

\section{Localization for the genus-\protect$ 0\protect $ GW invariants of toric
varieties}

In the setup described in the previous section,
Graber and Pandharipande's virtual localization formula applies:
if $ V $ is an equivariant vector bundle on the moduli space
$ \mathcal{M}_{0,m}(X_{\Sigma },A) $ then
\begin{equation}\label{eq:locgp}
\int\limits _{[\mathcal{M}_{0,m}(X_{\Sigma },A)]^{\text {virt}}}e^{T_{N}}(V)=
\sum _{\Gamma }
\frac{\iota _{*}e^{T_{N}}(V|_{\mathcal{M}_{\Gamma }})}{%
e^{T_{N}}(\mathcal{N}_{\Gamma }^{\text {virt}})}
\end{equation}
where $ \mathcal{N}_{\Gamma } $ is the so--called virtual normal bundle
to $ \mfix $ and $ e^{T_{N}} $ denotes the equivariant
Euler class. Since each $ \ev_{i}^{*}(Z_{j}) $ is equal to the (standard)
Euler class of an equivariant line bundle over $ \mathcal{M}_{0,m}(X_{\Sigma },A) $,
we obtain
\begin{equation}\label{eq:locGW}
\Phi _{0,m}^{X_{\Sigma },A}(Z^{l_{1}},\ldots ,Z^{l_{m}})=
\sum _{\Gamma }\frac{1}{\left| \mathbf{A}_{\Gamma }\right| }
\int\limits _{\mfix}
\frac{\prod\limits _{j=1}^{m}\prod\limits _{k=1}^{n}
\left( \omega _{k}^{\sigma (j)}\right) ^{l_{j,k}}}{%
e^{T_{N}}(\mathcal{N}_{\Gamma }^{\text {virt}})}
\end{equation}
where $ \omega _{k}^{\sigma (j)} $ is the following weight: Suppose
$ \sigma (j)=Z_{i_{1}}\cdots Z_{i_{d}} $
with $ i_{l}\neq i_{l'} $ whenever $ l\neq l' $. If $ k\notin \{i_{1},\ldots ,i_{d}\} $
then $ \omega _{k}^{\sigma (j)}:=0 $. Otherwise suppose (w.l.o.g.) $ i_{d}=k $ and
let $ \hat{\sigma }(j) $ be the unique vertex given by
$ \hat{\sigma }(j)=Z_{i_{1}}\cdots Z_{i_{d-1}}Z_{i_{d+1}} $
such that $ i_{d+1}\neq i_{l} $ for all $ l $. Then
$ \omega _{k}^{\sigma (j)}:=\omega_{\hat{\sigma }(j)}^{\sigma(j)} $.

By a careful analysis of the virtual normal bundle (see \cite[Theorem 7.2]{spielberg-phd})
we can compute its equivariant Euler class --- before we will give its formula here,
let us fix some notation. For an edge $e\in\Edge(\Gamma)$
define
\[
\Lambda _{\Gamma }(e):=\frac{(-1)^{h}h^{2h}}{(h!)^{2}(w^{\sigma }_{\sigma _{d}})^{2h}}
\prod _{j=1}^{d-1}
\frac{\prod\limits _{k=\lambda _{i_{j}}+1}^{-1}
\left( \omega _{\sigma _{j}}^{\sigma }-\frac{k}{h}\cdot \omega ^{\sigma }_{\sigma _{d}}\right)}{%
\prod\limits _{k=0}^{\lambda _{i_{j}}}
\left( \omega ^{\sigma }_{\sigma _{j}}-\frac{k}{h}\cdot \omega ^{\sigma }_{\sigma _{d}}\right)}.
\]
In this formula, we use the following notation: The edge $ e $ connects
the two fixed points $ \sigma  $ and $ \sigma _{d} $ with multiplicity
$ h $. The indices $ i_{j} $ and $ \hat{i}_{j} $ are chosen pair--wise
different such that $ \sigma =Z_{i_{1}}\cdots Z_{i_{d}} $ and
$ \sigma _{j}=Z_{i_{1}}\cdots Z_{\hat{i}_{j}}\cdots Z_{i_{d}} $ ($Z_{i_j}$ is
replaced by $Z_{\hat{i}_j}$).
The homology class of the edge $ e $ is given by
$ \lambda =(\lambda _{1},\ldots ,\lambda _{n}) $,
in particular $ \lambda _{i_{j}}=e\cdot Z_{i_{j}} $.

Furthermore let $ \Vertex_{t,s}(\Gamma ) $ be the set of vertices $v$ in the
graph $\Gamma$ with $ t $ outgoing edges and $ s $ marked points. For such
a vertex $v\in\Vertex_{t,s}(\Gamma)$ we define $\val(v):=t$, $\deg(v):=s+t$, and the class
\[ \omega(v):=
\left\{
\begin{array}{ll}
\omega_{F(v)}&\text{if $t=1$ and $s=0$}\\
1&\text{if $t=s=1$}\\
(\omega_{F_1(v)}+\omega_{F_2(v)})^{-1}&\text{if $t=2$ and $s=0$}\\
\left((\omega_{F_1(v)}-e_{F_1(v)})\cdots(\omega_{F_t(v)}-e_{F_t(v)})\right)^{-1}&
\text{if $t+s\geq 3$}.
\end{array}\right.
\]
Here the $e_{F}$ are Euler classes of universal cotangent lines to marked points
of $\mfix = \prod_{v\in\Vertex(\Gamma)} \DM[\deg(v)]$,
and the indices (from $1$ through $t$, for $\val(v)=t)$
refer to the different edges leaving the vertex $v$.

\begin{prop}[\cite{spielberg-phd}]
With this notation, the inverse of the equivariant Euler class of the virtual
normal bundle has the following expression:
\[
e^{T_N}(\mathcal{N}_\Gamma)^{-1} =
\left[ \prod_{t,s} \prod_{v\in\Vertex_{t,s}(\Gamma)}
\left( \omega_{\text{total}}^{\sigma(v)}\right)^{t-1} \cdot \omega(v)
\right] \cdot
\prod_{e\in\Edge(\Gamma)}\Lambda_\Gamma(e).
\]
\end{prop}

\begin{prop}[\cite{spielberg-GWcras,spielberg-GWpaper}]\label{prop:integralofnormal}
The integral over $\mfix$ of the invers of the equivariant Euler class of the virtual normal bundle
equals
\[
\int\limits _{\mfix}e^{T_{N}}(\mathcal{N}_{\Gamma })^{-1}=
\left[ \prod _{t,s}\prod _{v\in \Vertex_{t,s}(\Gamma )}
\frac{\left( \omega _{\text {total}}^{\sigma (v)}\right) ^{t-1}\cdot
\left( \omega _{\text {total}}^{F(v)}\right) ^{t+s-3}}{%
\prod _{i=1}^{t}\omega _{F_{i}(v)}}
\right]
\prod _{e\in \Edge(\Gamma )}\Lambda _{\Gamma }(e).
\]
\end{prop}

In particular, we see that if we want to generalize formula (\ref{eq:locGW})
to non--trivial classes $ \beta \in H^{*}(\DM,\mathbb Q) $,
the localization formula (\ref{eq:locgp})
tells us that it suffices to compute the equivariant Euler classes of the restrictions
of the equivariant bundles on $ \mathcal{M}_{0,m}(X_{\Sigma },A) $ representing
the class $ \beta  $, combine this class
with the equivariant Euler class of the virtual normal bundle, and integrate over
$\mfix$.

\section{Cotangent line bundles and their restrictions to the fixed point components}

In this section we will study how pull--backs of certain classes
$\beta\in H^*(\DM)$ localize to fixed--point
components $\mfix$.

Let $ \overline{\mathcal{C}}_{0,m}\longrightarrow \DM $
be the universal curve, and let
$ x_{i}:\DM\longrightarrow \overline{\mathcal{C}}_{0,m} $
be the marked point sections ($ i=1,\ldots ,m $). We will denote by
$ \mathbb L_{i}\longrightarrow \DM $
the $ i^{\text {th}} $ universal cotangent line, that is the pull--back by
$ x_{i} $ of the relative cotangent bundle
$ K_{\overline{\mathcal{C}}_{0,m}/\DM} $:
\[
\mathbb L_{i}:=x_{i}^{*}
\left( K_{\overline{\mathcal{C}}_{0,m}/\DM}\right) .
\]

For simplicity, we will restrict ourselves here to maximal sums of the line
bundles $\mathbb L_i$, i.e.\ to those of which the rank is equal to
$\dim \DM=m-3$. By Kontsevich's theorem
(\cite{kontsevich-Airyfunction}), we know that in this case
\begin{equation}\label{eq:kontsevich}
\int_{\DM}
e\left(\mathbb{L}_1^{\oplus d_1}\oplus\cdots\oplus\mathbb{L}_m^{\oplus d_m}\right)
= \frac{(m-3)!}{d_1!\cdots d_m!}
\end{equation}
that is the Euler class of this bundle is Poincar\'e dual to $(m-3)!/(d_1!\cdots d_m!)$
points; it is exactly this kind of classes $\beta$ we need in order to compute quantum
products of more than two factors (see equation (\ref{eq:multquantumproduct})). Note that
the $d_i$ fulfill  the equation $d_1+\cdots+d_m=m-3$.

\begin{lem}
The map $\pi:\mathcal{M}_{0,m}(X_\Sigma,A)\longrightarrow \DM$
forgetting the map to $X_\Sigma$ is equivariant with respect to the induced
$T_N$--action on $\mathcal{M}_{0,m}(X_\Sigma,A)$ and the trivial action on
$\DM$.
\end{lem}

\begin{proof}
Since the $T_N$--action on $\mathcal{M}_{0,m}(X_\Sigma,A)$ is induced from the
action on the image of the curve in $X_\Sigma$ (and which is discarded by the map
$\pi$), this is obvious.
\end{proof}

\begin{cor}\label{cor:trivialpullback}
The pull--back by $\pi$ of any bundle $E$ on $\DM$
is an equivariant bundle on $\mathcal{M}_{0,m}(X_\Sigma,A)$ with trivial fiber action.
\end{cor}

\begin{rem}
Corollary \ref{cor:trivialpullback}
implies in particular, that the equivariant and the non--equivariant
Euler classes of such pull--back bundles coincide; we will therefore use them
interchangeably in these cases.
\end{rem}

\begin{lem}
Let $E\longrightarrow \DM$ be a vector bundle of
$\rk E=m-3$. If $\pi(\mfix)\neq \DM$
then $e\left(\pi^* E|_{\mfix}\right)=0$.
\end{lem}

\begin{proof}
If $\pi(\mfix)\neq \DM$ then the codimension
of $\pi(\mfix)\subset \DM$ is at least
one. Therefore
\[
e\left(\pi^* E|_{\mfix}\right) = \pi^* e\left(E|_{\pi(\mfix)}\right)
= \pi^*(0) = 0,
\]
since $\rk E > \dim \pi(\mfix)$.
\end{proof}

The Lemma implies that if $e\left(\pi^* E|_{\mfix}\right)\neq 0$ for
a bundle $E$ with $\rk E=m-3$, the graph $\Gamma$ contains only one vertex
$v_\Gamma$ that corresponds to a stable component under the projection $\pi$
to $\DM$. In other words, if we fix $v_\Gamma$ as root
of the graph $\Gamma$, all its branches contain at most one marked point. We
will call such graphs $\Gamma$ {\em simple}.

\begin{thm}\label{thm:pullbackline}
Let $\Gamma$ be a simple $\mathcal{M}_{0,m}(X_\Sigma,A)$--graph, and let
$v_\Gamma$ be the unique stable vertex of $\Gamma$ and $\tilde{m}=\deg(v_\Gamma)$
be its degree. We will choose the indices of the marked points
$\tilde{x}_1,\ldots,\tilde{x}_{\tilde{m}}$ of $\DM[\tilde{m}]$
such that $\tilde{x}_1,\ldots,\tilde{x}_{m}$ are mapped by $\pi$ to
the marked points $x_1,\ldots,x_m$ in $\DM$, respectively.

Furthermore, let $\mathcal{L}_i\longrightarrow \DM[\tilde{m}]$ be the cotangent
lines to the marked points of $\DM[\tilde{m}]$. Then
\begin{multline}
e\left(\pi^*\left.\left(\mathbb{L}_1^{\oplus d_1}\oplus\cdots\oplus\mathbb{L}_m^{\oplus d_m}
\right)\right|_{\mfix}\right)
e\left( \mathcal{L}_{m+1}^{\oplus d_{m+1}} \oplus\cdots\oplus
\mathcal{L}_{\tilde{m}}^{\oplus d_{\tilde{m}}} \right)\\
=e \left(
\mathcal{L}_1^{\oplus d_1}\oplus\cdots\oplus \mathcal{L}_m^{\oplus d_m}
\oplus\cdots\oplus \mathcal{L}_{\tilde{m}}^{\oplus d_{\tilde{m}}} \right)
\end{multline}
whenever $d_1+\cdots+d_m=m-3$ and $d_{1}+\cdots+d_{\tilde{m}}=\tilde{m}-3$.
\end{thm}

\begin{proof}
First of all notice that $\pi$ factors through the projection map
$\mfix \longrightarrow \DM[\tilde{m}]$ to the factor corresponding
to the vertex $v_\Gamma$. Hence it suffices to consider the behaviour
of cotangent line bundles under $\pi_{\tilde{m},m}:\DM[\tilde{m}]
\longrightarrow \DM$, forgetting the last $(\tilde{m}-m)$ points. This
map $\pi_{\tilde{m},m}$ again factors into
\[
\pi_{\tilde{m},m}: \DM[\tilde{m}] \stackrel{\pi_{\tilde{m}}}{\longrightarrow}
\DM[\tilde{m}-1] \stackrel{\pi_{\tilde{m}-1}}{\longrightarrow} \cdots
\stackrel{\pi_{m+1}}{\longrightarrow} \DM.
\]
We will prove the statement for the map $\pi_{m+1}$ --- the Theorem then
follows by induction and Kontsevich's Theorem (equation (\ref{eq:kontsevich})).

So we want to show that (in the case $d_1+\cdots+d_m=m-3$)
\begin{multline}\label{eq:pullbacktoprove}
e\left( \pi_m^*\left(\mathcal{L}_1^{\oplus d_1}\oplus\cdots\oplus
\mathcal{L}_m^{\oplus d_m} \right) \oplus \mathcal{L}_{m+1}\right)\\
= e \left(
\mathcal{L}_1^{\oplus d_1}\oplus\cdots\oplus
\mathcal{L}_m^{\oplus d_m} \oplus \mathcal{L}_{m+1}\right).
\end{multline}
It is well--known (see e.g.\ \cite{harris-morrison}) that
\[ e(\mathcal{L}_i) = e(\pi_{m+1}^*(\mathcal{L}_i)) + D_i \]
where $D_i\cong \DM[3]\times\DM$ is the divisor in $\DM[m+1]$ where
the $\DM[3]$--bubble contains the marked points $x_i$ and $x_{m+1}$.
Now note that $\mathcal{L}_{m+1}|_{D_i}$ is constant for any $i=1,\ldots,m$,
hence $e(\mathcal{L}_{m+1})\cdot D_i=0$. This yields equation
(\ref{eq:pullbacktoprove}).
\end{proof}

\section{The formula for the Gromov--Witten invariants}

The next Corollary summarizes what we have shown so far:

\begin{cor}\label{cor:GWIfixed}
Let $E$ be the vector bundle $E=\mathbb{L}_1^{\oplus d_1}\oplus\cdots\oplus
\mathbb{L}_m^{\oplus d_m}$ on $\DM$ such that $d_1+\cdots+d_m=m-3$. Let
$\beta=e(E)$ be the Euler class of $E$. Then
\[
\Phi_{0,m}^{X_\Sigma,A}(\beta;Z^{l_1},\cdots,Z^{l_m}) =
\sum_{\Gamma \text{ simple}} \frac{1}{\left| \mathbf{A}_{\Gamma }\right| }
\int\limits _{\mfix}
\frac{e\left(\pi^* E|_{\mfix}\right)
\prod_{j,k=1}^{m,n}
\left( \omega _{k}^{\sigma (j)}\right) ^{l_{j,k}}}{%
e^{T_{N}}(\mathcal{N}_{\Gamma }^{\text {virt}})}.
\]
\end{cor}

We will now compute the integral over the fixed--point components to
obtain an explicit formula for these Gromov--Witten invariants.

\begin{thm}
Let $\Gamma$ be a simple $\mathcal{M}_{0,m}(X_\Sigma,A)$--graph, and let
$v_\Gamma$ be the unique stable vertex of $\Gamma$. Let $r:\Vertex(\Gamma)
\longrightarrow \mathbb N$ be the map defined by
\[
r(v):= \left\{ \begin{array}{ll}m-3&\quad \text{if $v=v_\Gamma$}\\
                                0&\quad \text{otherwise.}
               \end{array} \right.
\]
As before let $d_i$ be non--negative integers such that $d_1+\cdots+d_m=m-3$.
If we let $E=\mathbb{L}_1^{\oplus d_1}\oplus\cdots\oplus
\mathbb{L}_m^{\oplus d_m}$, then the following formula holds:
\begin{multline}\label{eq:newintegral}
\int_{\mfix} \frac{e^{T_N}\left(\pi^*(E)|_{\mfix} \right)}
{e^{T_N}\left(\mathcal{N}_\Gamma\right)}=
\frac{(\deg(v_\Gamma)-3)!}{d_1!\cdots d_m!(\deg(v_\Gamma)-m)!}\cdot
\prod_{e\in\Edge(\Gamma)} \Lambda_\Gamma(e)\cdot\\
\cdot\left[
\prod_{t,s} \prod_{v\in\Vertex_{t,s}(\Gamma)}
\frac{\left(\omega_{\text{total}}^{\sigma(v)}\right)^{t-1} \cdot
\left( \omega_{\text{total}}^{F(v)}\right)^{t+s-r(v)-3}}{%
\prod_{i=1}^t \omega_{F_i(v)}}
\right]_{\cdot}
\end{multline}

\end{thm}

\begin{proof}
By Kontsevich's Formula (\ref{eq:kontsevich}), it suffices to consider
$d_1=m-3$, $d_2=\ldots=d_m=0$. In this case, by Theorem \ref{thm:pullbackline},
the left hand side of (\ref{eq:newintegral}) equals
\[
\left[ \int_{\mfix} e\left(\mathcal{L}_{1,v_\Gamma}^{\oplus m-3}\right) \cdot
\prod_{t,s} \prod_{v\in\Vertex_{t,s}(\Gamma)}
\left(\omega_{\text{total}}^{\sigma(v)}\right)^{t-1} \cdot \omega(v)
\right]
\cdot
\prod_{e\in\Edge(\Gamma)} \Lambda_\Gamma(e).
\]
We therefore have to prove that
\begin{multline*}
\int_{\mfix}e\left(\mathcal{L}_{1,v_\Gamma}^{\oplus m-3}\right)
\cdot \prod_{v\in\Vertex(\Gamma)} \omega(v) =\\
=\frac{(\deg(v_\Gamma)-3)!}{(m-3)!(\deg(v_\Gamma)-m)!}
\prod_{t,s} \prod_{v\in\Vertex_{t,s}(\Gamma)}
\frac{\left(\omega_{\text{total}}^{F(v)}\right)^{t+s-r(v)-3}}{%
\prod_{i=1}^t \omega_{F_i(v)}}.
\end{multline*}
For the case that $v\in\Vertex_{t,s}(\Gamma)$ is different from $v_\Gamma$, we
have shown in the proof of Proposition \ref{prop:integralofnormal}
(see \cite{spielberg-GWpaper}) that
\[
\int_{\DM[\deg(v)]} \omega(v) = \frac{\left(
\omega_{\text{total}}^{F(v)}\right)^{t+s-3}}{\prod_{i=1}^t \omega_{F_i(v)}}.
\]
Hence we only have to consider the case when $v=v_\Gamma$; since it is very similar
to the previous case, we will only outline its proof. As in \cite{spielberg-GWpaper},
let $P_n(x_1,\ldots,x_k)=\sum_{\sum_i\tilde{d}_i=n} \prod_j x_j^{\tilde{d}_j}$. Let
$t:=\val(v_\Gamma)$, $\tilde{m}:=\deg(v_\Gamma)$ and $r:=m-3$. We will also
write $F_j$ instead of $F_j(v_\Gamma)$. Note that for the vertex $v_\Gamma$, we
always have $\tilde{m}-r-3\geq 0$.
Therefore
\begin{align*}
\int\limits_{\DM[\tilde{m}]} \omega(v_\Gamma) e\left(\mathcal{L}_1^{\oplus r}\right) &
= \int\limits_{\DM[\tilde{m}]} e\left(\mathcal{L}_1^{\oplus r}\right) \prod_{j=1}^t
\frac{1}{\omega_{F_j}-e_{F_j}}\\
&= \int\limits_{\DM[\tilde{m}]}e\left(\mathcal{L}_1^{\oplus r}\right) \prod_{j=1}^t
\frac{1}{\omega_{F_j}} \sum_{i=0}^\infty \left(\frac{e_{F_j}}{\omega_{F_j}}\right)^i\\
&\stackrel{(\ref{eq:kontsevich})}{=}
\prod_{j=1}^t \frac{1}{\omega_{F_j}}
 \sum_{\sum\tilde{d}_i=\tilde{m}-r-3} \frac{(\tilde{m}-3)!}{\tilde{d}_1!\cdots\tilde{d}_t!r!}
\left(\frac{1}{\omega_{F_1}}\right)^{\tilde{d}_1}\cdots
\left(\frac{1}{\omega_{F_t}}\right)^{\tilde{d}_t}\\
&=
\frac{(\tilde{m}-3)!}{r!(\tilde{m}-r-3)!}\cdot
\frac{\left(\omega_{\text{total}}^{F(v_\Gamma)}\right)^{\tilde{m}-r-3}}{%
\prod_{i=1}^t \omega_{F_i(v_\Gamma)}}
\end{align*}
which finishes the proof.
\end{proof}

\section{Example: The quantum cohomology of $\bsp$}

In this section we want to illustrate how useful the extension of the formula to the case
where $\beta=\pd(\pt)$ really is for the computation of the quantum cohomology ring
of a toric manifold.

In \cite{siebert-tian}, Siebert and Tian have shown that if the ordinary cohomolgy ring
of a symplectic manifold $X$ is given by
\[
H^*(X,\mathbb Q) = \mathbb{Q}[Z_1,\ldots, Z_k]/\langle R_1, \ldots, R_l\rangle
\]
where $R_1,\ldots,R_l$ are relations, then the quantum ring is given by
\[
QH^*(X,\mathbb Q) = \mathbb{Q}[Z_1,\ldots, Z_k]/\langle R^*_1, \ldots, R^*_l\rangle
\]
where $R_i^*$ are the relations $R_i$, but evaluated with respect to the quantum
product instead of the cup product.

Let us consider toric manifolds of the form $\bsp$. In \cite{costa-miro},
Costa and Miró--Roig have studied the three--point Gromov--Witten invariants
of these manifolds and announced that they will derive the quantum cohomology
ring of these manifolds in an upcoming paper. We have chosen the same example
to illustrate how the formula derived in this note can make computations
much easier.

We will recall some properties of $\bsp$ (for more details see
\cite{costa-miro}). Its cohomology ring is given by
\[
H^*(X_\Sigma, \mathbb Q)= \mathbb{Q}[Z_1,\ldots, Z_{r+2}]/\left\langle
{L_1, L_3, \ldots, L_{r+2}, P_1, P_2} \right\rangle.
\]
where the relations are given by
$L_1=Z_1-Z_2$, $L_3=Z_3-Z_{r+2}, \ldots, L_{r-1}=Z_{r-1}-Z_{r+2}$,
$L_r=Z_2+Z_r-Z_{r+2}$, $L_{r+1}=Z_2+Z_{r+1}-Z_{r+2}$,
$P_1=Z_1Z_2$, $P_2=Z_3\cdots Z_{r+2}$.

To see this, consider the fan whose one--dimensional cones are (with respect
to some basis $e_1,\ldots,e_r$ in the lattice ${\mathbb Z}^r$)
$v_1=e_1$, $v_2=-e_1+e_{r-1}+e_r$, $v_3=e_2,\ldots,v_{r+1}=e_r$ and
$v_{r+2}=-e_2-\cdots-e_r$, and
whose set of primitive collections is given by
\[
\mathfrak{P} = \left\{ \{v_1,v_2\},\{v_3,\ldots,v_{r+2}\}\right\}.
\]
In
\cite[Proposition 3.6]{costa-miro}, Costa and Miró--Roig
obtain (in what follows we will freely
use their notation)
\begin{equation}\label{eq:costamiro}
 Z_1 \star Z_1 = (Z_{r+2} \star Z_{r+2} - 2Z_1\star Z_{r+2}) \sum_{i\geq 1} q_1^i.
\end{equation}
Hence, we will only have to compute the quantum product $Z_3\star \cdots\star Z_{r+2}$
to get a presentation of the quantum cohomology ring. To do so, we will have to
compute the Gromov--Witten invariants of the form
\begin{equation}\label{eq:gwtocompute}
\Phi_{0,r+1}^{a\lambda_1+b\lambda_2}(\pd(\pt);Z_3,\ldots,Z_{r+2},\gamma).
\end{equation}

\begin{lem}
If $b\neq 1$, then all invariants of the form (\ref{eq:gwtocompute}) are zero.
\end{lem}

\begin{proof}
Suppose that for given $a$, $b$ and $\gamma$, the invariant (\ref{eq:gwtocompute})
is non--zero. Since $\langle c_1(X_\Sigma),a\lambda_1+b\lambda_2\rangle = 2rb$,
we must have $2r+\deg \gamma = 2rb+2r$. However, the real dimension of $X_\Sigma$
is equal to $2r$, hence $0\leq \deg\gamma \leq 2r$, and therefore $b\leq1$.

Now suppose $b$ is zero and therefore $\deg\gamma=0$ as well, that is $\gamma$ is
a multiple of the trivial class $1\in H^0(X_\Sigma)$. 
If $a=0$ as well, then the invariant is zero since it just equals the cup product
of the cohomology classes $Z_3,\ldots,Z_{r+2}$. So suppose $a>1$. Then
\begin{align*}
\lefteqn{\Phi^{a\lambda_1}_{0,r+1}(\pd(\pt);Z_3,\ldots,Z_{r+2},1)=}\\
&=\Phi^{a\lambda_1}_{0,r}(\pd(\pt);Z_3,\ldots,Z_{r-1},Z_{r+2},Z_r,Z_{r+1})\\
&=\sum_{a_1+\cdots+a_{r-2}=a \atop j_1,\ldots,j_{r-3}}
\Phi^{a_1\lambda_1}_{0,3}(Z_3,Z_{r+1},\gamma_{j_1})
\Phi^{a_2\lambda_1}_{0,3}(Z_4,\gamma_{j_1}^\vee,\gamma_{j_2})\cdots
\Phi^{a_{r-2}}_{0,3}(Z_{r+2},Z_r,\gamma_{j_{r-3}}^\vee),
\end{align*}
where the $\gamma_{j_i}$ run over a basis of $H^*(X_\Sigma,\mathbb Z)$. Since
$a>0$, at least one of $a_i$ has to be positive. On the other hand, we have
\[
\langle Z_3, \lambda_1\rangle=\ldots=\langle Z_{r-1},\lambda_1\rangle
=\langle Z_{r+2},\lambda_1\rangle = 0. 
\]
So as soon as $a_i>0$ the corresponding three--point invariant in the sum above
is zero. This proves the lemma.
\end{proof}

\begin{lem}
If $a>0$, then all invariants of the form (\ref{eq:gwtocompute}) are zero.
\end{lem}

\begin{proof}
Suppose that for given $a$, $b$ and $\gamma$, the invariant (\ref{eq:gwtocompute})
is non--zero. Then, by the previous Lemma, $b=1$. Hence we have to consider
homology classes $a\lambda_1+\lambda_2$, and $\gamma$ being a multiple
of the class of top degree, say $Z_1Z_3\cdots Z_{r}Z_{r+2}$.

The $1$--skeleton of the moment polytope of $X_\Sigma$ has the following
edges:
\begin{enumerate}
\item for $3\leq t < s \leq r+2$, edges between $\sigma_{1,t}$ and $\sigma_{1,s}$,
and between $\sigma_{2,t}$ and $\sigma_{2,s}$, all having homology class
$\lambda_2$;
\item for $3\leq t \leq r-1$ or $t=r+2$, an edge between $\sigma_{1,t}$ and
$\sigma_{2,t}$ of homology class $\lambda_1$;
\item for $t=r, r+1$, an edge between between $\sigma_{1,t}$ and
$\sigma_{2,t}$ of homology class $\lambda_1+\lambda_2$.
\end{enumerate}

The reader will now easily convince himself that there is no simple graph $\Gamma$
in this class
such that $Z_3,\ldots,Z_{r-1}$, $Z_{r+1}$, $Z_{r+1}$, $Z_{r+2}$ and $Z_1Z_3\cdots Z_{r},Z_{r+2}$
all have non--zero
equivariant Euler class on $\mfix$, unless $a=0$. Finally note that $Z_r=Z_{r+1}$
as cohomology classes, which finishes the proof.
\end{proof}

\begin{lem}
$\Phi^{\lambda_2}_{0,r+1}(\pd(\pt);Z_3,\ldots,Z_{r-1},Z_{r+1},Z_{r+1},Z_{r+2},
Z_1Z_3\cdots Z_rZ_{r+2})=1$.
\end{lem}

\begin{proof}
The only simple graph $\Gamma$ such that $Z_3,\ldots,Z_{r-1}$,
$Z_{r+1}$ $Z_{r+1}$, $Z_{r+2}$ and $Z_1Z_3\cdots Z_rZ_{r+2}$
all have non--zero equivariant Euler class on $\Gamma$ is the one
with one edge between $\sigma_{1,r+1}$ and
$\sigma_{1,r}$ where all but the last marked
point are at $\sigma_{1,r}$. Applying the formula derived in this note yields
the following for the Gromov--Witten invariant
$\Phi:=\Phi^{\lambda_2}_{0,r+1}(\pd(\pt);Z_3,\ldots,Z_{r-1},Z_{r+1},Z_{r+1},Z_{r+2},
Z_1Z_3\cdots Z_rZ_{r+2})$:
\begin{align*}
\Phi&=-
\frac{\omega_{\sigma_{2,r+1}}^{\sigma_{1,r+1}}\omega_{\sigma_{1,3}}^{\sigma_{1,r+1}}\cdots
\omega_{\sigma_{1,r}}^{\sigma_{1,r+1}}\omega_{\sigma_{1,r+2}}^{\sigma_{1,r+1}}
\omega_{\sigma_{1,3}}^{\sigma_{1,r}}\cdots\omega_{\sigma_{1,r-1}}^{\sigma_{1,r}}
\left(\omega_{\sigma_{1,r+1}}^{\sigma_{1,r}}\right)^2\omega_{\sigma_{1,r+2}}^{\sigma_{1,r}}}{%
\omega_{\sigma_{1,r+1}}^{\sigma_{1,r}}\left(\omega_{\sigma_{1,r+1}}^{\sigma_{1,r}}\right)^2
\omega_{\sigma_{1,r+2}}^{\sigma_{1,r}}
\left(\omega_{\sigma_{1,r+2}}^{\sigma_{1,r}}-\omega_{\sigma_{1,r}}^{\sigma_{1,r}}\right)
\omega_{\sigma_{1,r+2}}^{\sigma_{1,r+1}}}\cdot\\
&\phantom{===}\cdot\frac{1}{%
\omega_{\sigma_{1,3}}^{\sigma_{1,r}}
\left(\omega_{\sigma_{1,3}}^{\sigma_{1,r}}-\omega_{\sigma_{1,r+1}}^{\sigma_{1,r}}\right)
\cdots
\omega_{\sigma_{1,r-1}}^{\sigma_{1,r}}
\left(\omega_{\sigma_{1,r-1}}^{\sigma_{1,r}}-\omega_{\sigma_{1,r+1}}^{\sigma_{1,r}}\right)
}\\
&= \frac{\omega_{\sigma_{1,3}}^{\sigma_{1,r+1}}
\cdots \omega_{\sigma_{1,r-1}}^{\sigma_{1,r+1}}\omega_{\sigma_{1,r+2}}^{\sigma_{1,r+1}}}{%
\left(\omega_{\sigma_{1,3}}^{\sigma_{1,r}}-\omega_{\sigma_{1,r+1}}^{\sigma_{1,r}}\right)
\cdots
\left(\omega_{\sigma_{1,r-1}}^{\sigma_{1,r}}-\omega_{\sigma_{1,r+1}}^{\sigma_{1,r}}\right)
\left(\omega_{\sigma_{1,r+2}}^{\sigma_{1,r}}-\omega_{\sigma_{1,r+1}}^{\sigma_{1,r}}\right)}\\
&=1
\end{align*}
where we have used $\omega_{\sigma_{1,i}}^{\sigma_{1,r+1}}=
\omega_{\sigma_{1,i}}^{\sigma_{1,r}}-\omega_{\sigma_{1,r+1}}^{\sigma_{1,r}}$ for
$i=3,\ldots,r-1,r+2$ by a similar argument as in \cite[Lemma 2.1]{costa-miro},
that is by applying \cite[Lemma 6.7]{spielberg-phd}.
\end{proof}

\begin{cor}
The quantum cohomology ring of $X_\Sigma=\bsp$ (for $r\geq 3$) is given by
\[
QH^*(X_\Sigma, \mathbb Z) = {\mathbb Z}[Z_1,\ldots, Z_{r+2},q_1,q_2] /
\left\langle L_1,L_3,\ldots,L_{r-1},R_1^*,R_2^*\right\rangle,
\]
where the relations $R_1^\star$ and $R_2^\star$ are given by
\begin{align*}
R_1^\star&=Z_1 \star Z_2-Z_{r}\star Z_{r+1} q_1\\
R_2^\star&=Z_3\star\cdots\star Z_{r+2} - q_2,
\end{align*}
that is this quantum ring coincides with Batyrev's ring in
\cite{batyrev-quantumcohomologytoric}.
\end{cor}

\begin{proof}
The relation $Z_3\star \cdots\star Z_{r+2}=q_2$ follows directly from
the previous three Lemmata. Now let us consider the other multiplicative
relation:
\begin{align*}
Z_1 \star Z_2 &= Z_1 \star Z_1 \quad \text{since $Z_1=Z_2$}\\
&= \sum_{i\geq 1} (Z_{r+2} \star Z_{r+2} - 2 Z_1\star Z_{r+2}) q_1^i \quad
\text{by equation (\ref{eq:costamiro})}\\
&= (Z_{r} \star Z_{r+1} - Z_1 \star Z_2) \sum_{i\geq1} q_1^i\quad \text{by the linear relations}.
\end{align*}
Hence we obtain
\[ Z_1 \star Z_2 \sum_{i\geq 0} q_1^i = Z_{r} \star Z_{r+1} \sum_{i\geq 1} q_1^i,
\]
and by comparing the coefficients of the $q_1^i$, the relation follows.
\end{proof}

\begin{rem}
Note that although we derive the same presentation for the quantum
cohomology ring as the one stated by Batyrev in
\cite{batyrev-quantumcohomologytoric},
the Gromov--Witten invariants that enter as structure constants
into the computation of the quantum products are different from
those numbers considered by Batyrev.
\end{rem}

\section*{Acknowledgements}

I am grateful to Laura Costa and Rosa M. Miró--Roig for having sent and explained their paper
\cite{costa-miro} to me.
I also want to thank the Centro de Ana\-líse Matemática, Geometria e Sistemas
Dinâ\-mi\-cos at
Instituto Superior Técnico, Lisboa, for its hospitality.

\bibliographystyle{alpha}
\bibliography{mathbib}

\begin{thebibliography}{CMR01}

\bibitem[Bat93]{batyrev-quantumcohomologytoric}
Victor~V. Batyrev.
\newblock Quantum cohomology rings of toric manifolds.
\newblock {\em Ast\'erisque}, 218:9--34, 1993.

\bibitem[Beh97]{behrend-gwinvariants}
Kai Behrend.
\newblock {Gromov--Witten} invariants in algebraic geometry.
\newblock {\em Invent. Math.}, 127:601--627, 1997.

\bibitem[CMR01]{costa-miro}
Laura Costa and Rosa Miró-Roig.
\newblock {GW}-invariants and quantum products with infinitely many quantum
  corrections.
\newblock Preprint, to appear in forum mathematicum, Universitat de Barcelona,
  June 2001.

\bibitem[Ful93]{fulton-toricvarieties}
William Fulton.
\newblock {\em Introduction to toric varieties}.
\newblock Number 131 in Annals of Mathematics Studies. Princeton Univ.~Press,
  1993.

\bibitem[GP99]{graber-pandharipande}
T.~Graber and R.~Pandharipande.
\newblock Localization of virtual classes.
\newblock {\em Invent. math.}, 135(2):487--518, 1999.

\bibitem[HM98]{harris-morrison}
Joe Harris and Ian Morrison.
\newblock {\em Moduli of curves}.
\newblock Number 187 in Graduate Texts in Math. Springer Verlag, 1998.

\bibitem[Kon92]{kontsevich-Airyfunction}
Maxim Kontsevich.
\newblock Intersection theory on the moduli space of curves and the matrix
  {Airy} function.
\newblock {\em Comm.\ Math.\ Phys.}, 147:1--23, 1992.

\bibitem[Oda88]{oda-convexbodies}
Tadao Oda.
\newblock {\em Convex bodies and algebraic geometry}.
\newblock Number~15 in Ergebnisse der Mathematik und ihrer Grenzgebiete,
  3.~Folge. Springer Verlag, 1988.

\bibitem[RT95]{ruan-tian-mathematicaltheory}
Yongbin Ruan and Gang Tian.
\newblock A mathematical theory of quantum cohomology.
\newblock {\em J.\ Diff.\ Geom.}, 42(2):259--367, 1995.

\bibitem[Spi99a]{spielberg-phd}
Holger Spielberg.
\newblock A formula for the {Gromov--Witten} invariants of toric varieties.
\newblock PhD thesis/Preprint 1999/11, IRMA, Universit\'e Louis Pasteur,
  Strasbourg, February 1999.

\bibitem[Spi99b]{spielberg-GWcras}
Holger Spielberg.
\newblock The {Gromov--Witten} invariants of symplectic toric manifolds, and
  their quantum cohomology ring.
\newblock {\em {C.~R.~Acad.~Sci.~Paris}, S\'erie I}, 329(8):699--704, 1999.

\bibitem[Spi00]{spielberg-GWpaper}
Holger Spielberg.
\newblock The {Gromov--Witten} invariants of symplectic toric manifolds.
\newblock Preprint math.ag/0006156 at arxiv.org, June 2000.

\bibitem[ST97]{siebert-tian}
Bernd Siebert and Gang Tian.
\newblock On quantum cohomology rings of {Fano} manifolds and a formula of
  {Vafa} and {Intriligator}.
\newblock {\em Asian J. Math.}, 1(4):679--695, 1997.

\bibitem[Wit91]{witten-2dimgravity}
Edward Witten.
\newblock Two--dimensional gravity and intersection theory on moduli spaces.
\newblock {\em Surveys Diff.\ Geom.}, 1:243--310, 1991.

\end{thebibliography}

\end{document}